\documentclass{aptpub}

\authornames{V. ANTIPOV  AND Yu. KABANOV}

\shorttitle{Ruin Probabilities with Investments in  Random Enviroments: Smoothness} 

\numberwithin{equation}{section}

\usepackage{amscd,amsfonts,amssymb,amsmath,latexsym,array,hhline,xcolor,graphicx}

\newcommand\cE{{\cal E}}

\newcommand\cF{{\cal F}}

\newcommand\cI{{\cal I}}

\newcommand\cL{{\cal L}}

\newcommand\e{{\varepsilon}}

\def\E{{\bf E}}
\def\P{{\bf P}}
\def\Chi{{\bf 1}}

\def\text#1{\hbox{#1}}

\def\E{{\bf E}}
\def\P{{\bf P}}

\def\Chi{{\bf 1}}

\def\build #1_#2{\mathrel{\mathop{\kern 0pt #1}\limits_{#2}}} 

\def\R{{\mathbb R}}

\newcommand\fdem{$\Box$}

\begin{document}


\title{Ruin Probabilities with Investments in  Random Environment: Smoothness}

\authorone[Lomonosov Moscow State University and ``Vega" Institute]{Viktor Antipov} 
\emailone{stayaptichek@gmail.com}

\authortwo[Lomonosov Moscow State University and Universit\'e de Franche-Comt\'e]{Yuri Kabanov}

\addresstwo{16 Route de Gray,  25030 Besan\c{c}on cedex, France}
\emailtwo{ykabanov@univ-fcomte.fr} 


\begin{abstract}
The paper deals with the ruin problem  of an insurance company 
 investing  its capital reserve in a risky asset with the price dynamics given by a conditional geometric Brownian motion whose parameters 
 depend on a Markov process describing  a random  variations in the economic and financial environments.  
 We prove  smoothness of the ruin probability as  a function of the  initial capital and  obtain for it an
integro-differential equation.   

\end{abstract}

\smallskip


\smallskip

\newcommand\beq{\begin{equation}}
\newcommand\eeq{\end{equation}}
\newcommand\bea{\begin{eqnarray}}
\newcommand\eea{\end{eqnarray}}
\newcommand\bean{\begin{eqnarray*}}
\newcommand\eean{\end{eqnarray*}}

 \section{Introduction}  

Ruin models with risky investments (in other words, with stochastic interest rates) is one of the most active fields of the present day ruin theory. 
These models combine classical frameworks  with models of asset price dynamics developed in mathematical finance. The goal of numerous studies is to provide an information about asymptotic of ruin probabilities. 
For this there are several approaches. One of them is based on the integro-differential equations for the ruin probabilities. An inspection of the  literature reveals that in many cases  these equations are derived  assuming the smoothness, see, e.g., \cite{Paul-93}, \cite{RP}.  To our mind, the smoothness of ruin probabilities  is a rather delicate property  requiring a special (and rather involved)  study but there are very few papers on it, see \cite{KP2016} and \cite{KPukh}.    

In the present note we consider the price dynamics   suggested  by Di Masi, Kabanov, and Runggaldier,  \cite{DKR}, where the coefficients of a (conditional) geometric Brownian motion depend on a Markov process with finite number of states.  Such a setting (sometimes referred to as stochastic volatility model, regime switching, or hidden Markov chain)  reflects the random dynamics of the economic environment and  seems   to be adequate for a context with long term contracts typical in insurance. It was already considered in the actuarial  literature, see, e.g., the papers
\cite{EK}, \cite{KP2022} where the analysis was based using implicit renewal theory and  the paper \cite{RP} where a  second order integro-differential equation for the ruin probability was derived using intuitive arguments.   

Our main result is the theorem asserting that if the distribution of jumps has a density with two continuous and integrable derivatives then the ruin probability is two times continuously differentiable.

\section{The model}
Let   $(\Omega,\cF,{\bf F}=(\cF_t)_{t\in \R_+},\P)$  be a stochastic basis where are given  three independent stochastic processes, $W$, $P$, and $\theta^i$ generating the filtration ${\bf F}$:

1. A standard Wiener process $W=(W_t)$.

2.  A compound Poisson process $P=(P_t)$ with drift:   
$$
P_t=ct+\sum_{k=1}^{N_t}\xi_k
$$
where $c$ is a constant, $N=(N_t)$ is a Poisson process with intensity $\alpha>0$, and $(\xi_k)$ is an i.i.d. sequence independent of $N$  with 
the distribution $F_\xi=F_\xi(dx)$.  Sometimes, it is more convenient   the alternative description using  the notation of  stochastic calculus, namely,
\beq
 \label{r}
P_t=ct+x*\pi_t = ct  + \int_0^t\int x  \pi(ds,dx)  
\eeq
where $\pi=\pi(dt,dx)$ is the jump measure of $P$, that is a  Poisson random measure  $\R_+\times \R$ with the deterministic  compensator 
$\tilde \pi=\tilde \pi(dt,dx)=\Pi(dx)dt$ (coinciding  with the mean of $\pi$) and the L\'evy measure $\Pi(dx)=\alpha F_\xi(dx)$. We denote  by $T_n$the consecutive  jumps of the process $N_t=\pi([0,T]\times \R)$ with the usual convention  $T_0=0$. 

3.  A piecewise constant right-continuous  Markov process $\theta^i=(\theta^i_t)$  with values in  the finite set  $E:=\{0,1,...,K-1\}$, the initial value $\theta^i_0=i$, and the 
transition intensity matrix $(\lambda_{jk})$. We  assume that all states are communicating.   

 Recall that $\lambda_{jj}=-\sum_{k\neq j}\lambda_{jk}$ and $\lambda_j:=-\lambda_{jj}>0$ for each $j$. 
We denote $\tau^i_n$   the consecutive jumps  of $\theta^i$ with the convention  $\tau^i_0=0$ and  
 consider also the imbedded Markov chain $\vartheta^i_n:=\theta^i_{\tau^i_n}$ with 
transition probabilities $P_{kl}=\lambda_{kl}/\lambda_{k}$, $k\neq l$, and $P_{kk}=0$.

The conditional distribution  of the length of interval $\tau^i_{k+1}-\tau^i_{k}$   with respect to the $\sigma$-algebra $\cF_{\tau^i_k}$
is exponential with parameter $\lambda_{\vartheta^i_k}$. 
 
We consider  the random integer-valued measure 
$p^i:=\sum_n\e_{(\tau^ i_n,\vartheta_n^i)}$ where $\e_x$ is the unit mass at $x$.  It can be also   defined by  the counting processes   
$$
p^i([0,t],k):=\sum_{\tau^ i_n\le t}\Chi_{\{\vartheta_n^i=k\}}, \qquad k\in E.
$$
The compensator $\tilde p^i$ (dual predictable projection) of the random measure $p^i$ is given by the formula
\beq
\label{compensator}
\tilde p^i(dt,k,)=\Chi_{\{ \theta^i_{t}\neq k\}}\lambda_{\theta^i_{t}k}dt.
\eeq

To alleviate formulae we shall skip the superscript $i$ (and also u) when it will not lead to an ambiguity and use for brevity the standard ``dot" notation of the semimartingale stochastic calculus for integrals, with   $L=(L_t)$ standing for the (deterministic) process with $L_t:=t$.  

With these conventions the Markov modulated price process $S$ can be given as the solution of the linear integral 
equation 
$$
S=1+S(a_{\theta}\cdot L+\sigma_{\theta}\cdot W). 
$$
In the traditional symbolical  form of the Ito calculus this is written as    
$$
dS_t=S_t (a_{\theta_t} dt+\sigma_{\theta_t}dW_t), \qquad S_0=1,    
$$ 
where $a_k\in \R$, $\sigma_k>0$, $k=0,...,K-1$.   

Let $\kappa_k:=a_k-(1/2)\sigma^2_k$. 
The solution  of (homogeneous) linear equation can be represented 
by the formula $S=\cE(R)$ where $\cE$ is the stochastic exponential and $R:=a_{\theta}\cdot L+\sigma_{\theta} \cdot W$ is the {\it relative price process}, or, in a more explicit form, by the formula $S=e^V$ where 
$$
V:=\kappa_\theta\cdot L + \sigma_\theta\cdot W
$$ is the {\it logprice process}. 

The  process $X=X^{u,i}$ is defined by the formula  
\beq
 \label{risk}
 X=u+  X  \cdot R+ P   
 \eeq
 where $u>0$.

In the actuarial context  the process  $X=X^{u,i}$ is interpreted  as the capital  reserve 
 of an insurance company   fully invested  in a risky asset whose price 
$S^i$ is  a conditional geometric Brownian motion given the Markov process $\theta^i$ describing the financial environment.  

The process $P$ describes the business activity.  There are two basic models: of the non-life insurance ($c<0$, $F_\xi((0,\infty))=1$) and of the  annuity payments ($c>0$, $F_\xi((-\infty,0))=1$.  In the literature one can find also a mixed model where 
$F_\xi$ charges both half-axes.  

We consider here the setting where only the coefficients of the price process (i.e. the stochastic interest rate) depend on $\theta$.  The extension   to the case where  parameters of the business process also depends on $\theta$  is rather straightforward and  has no effect on our main results.    

We  assume that 
$P$ is not increasing  (otherwise the ruin probability is zero).  
 
\smallskip 
 Let $\tau^{u,i}:=\inf \{t>0\colon X^{u,i}_t\le 0\}$ be the instant of ruin corresponding to the initial capital $u$ and the initial regime $i$. 
 Then 
 $\Psi_i(u):=\P[\tau^{u,i}<\infty]$ is the ruin probability  
 and $\Phi_i(u):=1-\Psi_i(u)$
 is the survival probability.

\smallskip 

Our main result is a sufficient condition on  the smoothness of survival (and ruin) probabilities as functions of the initial capital.  
\begin{theorem}
 Let $\xi>0$ be a random variable with a density $f$  which is two times  continuously  differentiable on $(0,\infty)$  and such that   $f',f''\in L^1(\R_{+})$.
Then the functions $\Phi^i(u)$ are  two times  continuously differentiable for $u>0$.
\end{theorem}

The proof is given in Section \ref{smooth}. 

\section{Smoothness of the survival probability}
\label{smooth}

To use the theory of Markov processes we should consider not a single  process $\theta$ but a 
family 
of processes $\theta^i$ with initial value $\theta^i_0=i$, $i\in E$. Automatically,  other related processes also will  depend on $i$, in many cases omitted in the notations to alleviate formulae.   

The basic idea to prove the smoothness of the survival probability   is to use an integral representation.     
Define the continuous process $Y=Y^{u,i}$ with 
\beq
\label{Yuit}
    Y^{u, i}_t:= e^{V^i_t} \left( u + c\int_0^t e^{-V^i_s} ds\right), 
\eeq
coinciding with $X^{u, i}$ on $[0, T_1)$.

 We consider  the case of non-life insurance, i.e. with $c>0$. 

By virtue of the strong Markov property of  $(X^{u, i},\theta^i)$
$$
\Phi_{i}(u)  =\sum_{j\in E} \E \left[\Phi_j(X^{u, i}_{ T_1}) \Chi_{\{\theta^i_{T_1}=j\}}\right]=\sum_{j\in E} \E \left[\Phi_j(X^{u, i}_{ T_1-}+\Delta X^{u, i}_{T_1}) \Chi_{\{\theta^i_{T_1}=j\}}\right].
$$
On the interval $[0,T_1)$ the process $X^{u, i}$ coincides with the continuous increasing process $Y^{u, i}$ and 
$\Delta X^{u, i}=\xi_1<0$. Putting $\alpha(dt):= \alpha e^{-\alpha t}dt$ we get that  
\bean
\Phi_i(u)&=&\sum_{j\in E}\E\left[ \Chi_{\{\theta^i_{ T_1}=j\}}\Phi_j(Y^{u,i}_{T_1}+\xi_1)\right]
=\sum_{j\in E}(\Upsilon^{1i}_j(u)+  \Upsilon_j^{2i}(u)),
\eean
where
\bean
\Upsilon^{1i}_j(u)&:= &\int_0^2\E \left[\Chi_{\{\theta^i_{ t}=j\}} \Phi_j(Y^{u,i}_t+\xi_1)\right]\alpha(dt),\\
\Upsilon^{2i}_j(u)&:=& \int_2^\infty\E \left[\Chi_{\{\theta^i_{ t}=j\}} \Phi_j(Y^{u,i}_t+\xi_1)\right]\alpha(dt).
\eean
Clearly, 
\bean
\Upsilon^{1i}_j(u)&:= &\int_0^2\int_0^\infty  \int_{D}  \E\left[  \Phi_j(Y_t^{u}(z)-w)\right]
\Chi_{\{z_{ t}=j\}}m^ i(dz)F_{|\xi|}(dw)\alpha(dt),\\
\Upsilon^{2i}_j(u)&:=&\int_2^\infty \int_{D}  \E\left[\tilde \Phi_j(Y_t^{u}(z))\right]\Chi_{\{z_{ t}=j\}}m^ i(dz)(dw)\alpha(dt),
\eean
where $m^i(dz)$ is the distribution of the process $\theta^i$ in the Skorohod space  $D$ of the  c\`adl\`ag functions with a generic point $z=(z_s)$ and 
 $ \tilde \Phi_j(y):= \E\Phi_j(y-|\xi_1|)$. Recall that we consider the case where $\xi_1=-|\xi_1|$.

\medskip 
\subsection{Smoothness of $\Upsilon^2$.}

In this case the distribution $F_\xi$ is not involved (and this is the reason why we introduce the function $\tilde \Phi_j$). 

\begin{lemma}
Let $G:\R\to [0,1]$ be a Borel function such that $G(y)=0$ for $y\le 0$. Let  $J(u,t,z):=\E\left[ G(Y_t^{u}(z))\right]$, $t\ge 2$. Then 
$J(.,t,z)\in C^{\infty}((0,\infty))$ and $$
\sup_z(\partial^n/\partial u^n)J(u,t,z)\le C(u).
$$ 
\end{lemma}
{\sl Proof.} To simplify formulae we write here ``local" notations  $V$, $\sigma$, etc. instead of $V(z)$, $\sigma(z)$, etc. 
Using the representation  
$$
Y_t^u= e^{V_t-V_1}Y_1^u+c \int_1^t e^{V_t-V_s} ds, \qquad t\ge 1,  
$$
and  the independence of the process $(V_s-V_1)_{s\ge 1}$  and the random variable   $Y_1^u$ we get that 
$$ 
\E [G(Y^u_t)|Y^u_1]=G(t,Y_1^u),
$$
where
$$  
G(t,y) := \E \left[G\left (e^{ V_{t}- V_{1}}y+
c\int_{1}^{t}
e^{ V_{t} -V_{s}}d s\right)\right]. 
 $$ 
Substituting the expression for $Y^u_1$ given by (\ref{Yuit}) 
and  using abbreviated notations for the integrals
$$
\kappa\cdot \Lambda_t:=\int_0^t\kappa_sds, \quad \sigma\cdot W_t=\int_0^t \sigma_sdW_s, 
\quad 
M_t:=ce^{-\kappa\cdot \Lambda }\cdot \Lambda_t=c\int_0^te^{-\int_0^s \kappa_rdr}ds, 
$$
we have
\bean
J(u,t)&:=&\E\left[ G(Y_t^{u} ) \right]=\E\left[ G(t,Y_1^{u} ) \right]= 
\E\left[ G\left(t,e^{\kappa\cdot \Lambda_1+\sigma\cdot W_{1}}\big (u+e^{-\sigma\cdot W}\cdot M_1\big)\right)\right] \\
&=&\E\left[ \tilde G\left(t,e^{\sigma\cdot W_{1}}\big (u+e^{-\sigma\cdot W}\cdot M_1\big)\right)\right],
\eean
where 
$  
\tilde G(t,y) := G(t,e^{\kappa\cdot \Lambda_1}y) $.

The next step is a reduction to the case of constant coefficients.  The arguments are simple.
The function 
$A= \sigma^2\cdot \Lambda$ is strictly increasing. Let us denote $A^{-1}$ its inverse. 
 The law of the continuous process with independent increments $(\sigma \cdot W_s)_{s\ge 0}$
coincides with the law of the process $(W_{A_s})_{s\ge 0}$. 

Let $m(s):=  c e^{-\kappa\cdot \Lambda_s}$, $\tilde m(r) := m(A^{-1}_r)/\sigma^
2_{A^{-1}_{r}}$. Changing the variable we get that 
 $$
  \int_{0}^1 e^{-W_{A_s}}dM_s=  \int_{0}^1 e^{-W_{A_s}}m(s)ds=\int_{0}^{A_1}e^{-W_{r}}\tilde m(r)dr. 
 $$
It follows that 
$$
J(u,t)=\E\left[ \tilde G\left(t,e^{W_{A_1}}\big (u+e^{-W}\tilde m\cdot \Lambda_{A_1}\big)\right)\right]. 
$$
Let $T_1(s):=A^{-1}_{A_1s}$ and $\bar \sigma:=\sqrt {A_1}$. Then $(\bar \sigma W_{t/\bar \sigma ^2})_{t\ge 0}$ is a Wiener process and  
 we get that 
$$
J(u,t)=\E\left[ \tilde G\left(t,e^{\bar \sigma W_1}\big (u+e^{-\bar\sigma  W}\bar m\cdot \Lambda_{1}\big)\right)\right]
$$
where $\bar m(s):=\tilde m(A_1s)A_1 = m(T_1(s))A_1/\sigma^2_{T_1(s)}$.

Using the conditioning with respect to $W_1=x$  and taking into account that the conditional law of $(W_t)_{t\in[0,1]}$ in this case is that of the Brownian bridge on $[0,1]$ ending at the level $x$, that is the same as the process
$W_s+s(x-W_1)$, $s\in [0,1]$,  
 we get that 
 $$
 J(u,t)=\int \E [H(t,x,u+\eta_x)
]\varphi_{0,1}(x)dx
 $$
where $H(t,x,y)=\tilde G(t,e^{\bar \sigma x}y)$ is a function taking values in $[0,1]$,  
\bea
\eta_x&:=&\int_0^1 e^{- \bar\sigma (W_s-sW_1)} m(s,x)ds, \\
 m(s,x)&:=&e^{-\bar\sigma sx}\bar m(s)=e^{-\bar\sigma sx} m(T_1(s))  \bar \sigma^2 /\sigma^2_{T_1(s)}.
\eea

Let us check that $\eta_x$ has a density. To this end we consider on $[0,1]$  the Gaussian process $D_s:=W_s-sW_1-\gamma_sW_{1/2}$ where  
$\gamma_s:=(2s)\wedge 1-s$. 
Note that 
$$
\E [D_sW_{1/2}]=\E [(W_s-sW_1-\gamma_sW_{1/2})W_{1/2}]=s\wedge 1/2-s/2-\gamma_s/2=0
$$
and, hence,  $D$ and $W_{1/2}$ are independent. 
It follows that  for any bounded Borel function $g$  
$$
\E [g(\eta_x)
]=\E \left[\int g(B_x(v))\varphi_{0,1/2}(v)dv\right]
$$
where the strictly positive  $C^\infty$-function 
$$
v\mapsto B_x(v):=\int_0^1e^{-\bar \sigma (D_s+ \gamma_s v)} m(s,x)ds
$$
(depending of $x$ and also of $\omega$, omitted as usual) is  strictly decreasing, $B_x(-\infty)=\infty$ and  $B_x(\infty)=0 $.  The inverse function $y\mapsto B_x^{-1}(y)$ is also strictly decreasing $B_x^{-1}(0+)=\infty$  and $B_x^{-1}(\infty)=  - \infty$. 

 After the change of variable we get that 
$$
\E [g(\eta_x)]= \int_0^\infty g(y)\E\left[\frac { \varphi_{0,1/2}(B^ {-1}_x(y))}{K_x(B^ {-1}_x(y))}\right]dy
$$
where 
\beq
\label{Kxv}
K_x(v)=B_x'(v)= -\bar \sigma \int_0^1 \ \gamma_se^{- \bar \sigma (D_s+\gamma_s v)}m(s,x)ds.
\eeq
Thus, $\eta_x$ has the density 
$$
\rho_x(y)=\E\left[\frac { \varphi_{0,1/2}(B^ {-1}_x(y))}{K_x(B^ {-1}_x(y))}\right], \qquad y>0.  
$$
Changing the variable we get that 
$$
  J(u,t)=\int \E[ H(t,x,u+\eta_x)]\varphi_{0,1}(x)dx= \int_{\mathbb{R}^2} H(t,x,y)\rho_x(y - u)\varphi_{0,1}(x)dydx.  
$$ 
It remains to show that $\rho_x(.)$ belongs to $C^\infty$  and its derivatives has at most exponential growth. 

Put 
$$Q^{(0)} (x, v)  := \frac{\varphi_{0, 1/2} (v)}{K_x(v)}, \quad Q^{(n)} (x, v)  := - \frac{\partial Q^{(n-1)}(x, v)}{\partial v}\frac 1{K_x(v)}, \quad n \geq 1.$$
Then 
\bean
\frac{\partial}{\partial y} Q^{(0)} (x, B^{-1}_{x}(y) ) &=& \frac{\partial Q^{(0)}(x,B^{-1}_{x}(y))}{\partial v} \frac{\partial B^{-1}_{x}(y)}{\partial y}= - \frac{\partial Q^{(0)}(x, B^{-1}_{x}(y))}{\partial v}\frac 1{K_x(B^{-1}_{x}(y))}\\
& =& Q^{(1)} (x, B^{-1}_{x}(y) )
\eean
and, similarly,
$$
\frac{\partial}{\partial y} Q^{(n-1)} (x, B^{-1}_{x}(y) ) = - \frac{\partial Q^{(n-1)} (x, B^{-1}_{x}(y))}{\partial v}\frac 1{K_x(B^{-1}_{x}(y))} = Q^{(n)} (x, B^{-1}_{x}(y) ).
$$
It follows that 
$$
\frac{\partial^n}{\partial y^n} Q^{(0)} (x, B^{-1}_{x}(y) )= Q^{(n)} (x, B^{-1}_{x}(y) ). 
$$
It is easily seen that 
$$
Q^{(n)}(x, v) = \frac{\varphi_{0, 1/2} (v)}{K^{n+1}_{x} (v)} \sum_{k=0}^{n} P_k(v) R_{n-k} (x, v),
$$
where $P_k(v)$ is a polynomial of order $k$ and $R_{n-k}(x, v)$ is a linear combination of products of derivatives of $K_x(v)$ in variable $v$. 

Recall that 
$$
m(s,x) =e^{-\bar\sigma sx}  m(T_1(s))\bar \sigma^2/\sigma^2_{T_1(s)}=
 c\,  e^{-\bar \sigma sx - \kappa \cdot \Lambda_{T_1(s)}}\bar \sigma^2/\sigma^2_{T_1(s)}.
$$
It follows that for all $s\in [0,1]$ and $x\in \R$ we have the bounds
$$
c_me^{-\bar \sigma |x|}\le  m(s,x)\le C_me^{\bar \sigma |x|}
$$
where $c_m,C_m>0$ are  constants depending only on bounds on $|a|$ and $\sigma$. 
Let $\sigma^*:=\sup_s \sigma_s$,  $\kappa^*:=\sup_s |\kappa_s|$, and    $W_1^{*} := \sup_{s \in [0,1]} |W_s|$.
Note that $\gamma\le 1/2$ and 
$$
\gamma^n\cdot \Lambda_1=2^{-(n+1)}/(n+1), \quad n\ge 1.
$$   
Clearly,  
$$
e^{-\sigma^*(3W^*_1 +v)}\le 
e^{- \bar \sigma (D_s+\gamma_s v)}\le e^{\sigma^*(3W^*_1 +v)}, 
$$
Combining this with the above bounds we easily get from  (\ref{Kxv}) that 
$$
 \frac 14 \sigma_* c_m e^{-\sigma^*(|x|+|v|+3 W_1^*)}\le
  |K_x(v)|
 \le   \frac 14 \sigma^* C_m e^{\sigma^*(|x|+|v|+3 W_1^{*})},
$$
$$
\frac { (\sigma_*/2)^{n+1}}{n+2}  c_m e^{-\sigma^*(|x|+|v|+ 3 W_1^*)}\le 
 \left |\frac{\partial^n}{\partial v^n}K_x(v)\right|
 \le  \frac {(\sigma^*/2)^{n+1}}{n+1}  C_m e^{\sigma^*(|x|+|v|+3  W_1^{*} )}. 
$$

It follows that there exists a constant $C_n > 0$ such that 
\beq
\label{Qn}
|Q^{(n)}(x, v)| \leq C_n e^{C_n |x|} (1 + v^n) e^{-v^2} e^{C_n W_1^{*}}.
\eeq
Therefore, 
$$ 
\E \left[\sup_y |Q^{(n)} (x, B^{-1}_{x}(y) )|\right] \le  \E  \left[  \sup_{v \in \mathbb{R}}\left |  Q^{(n)}(x, v)\right |\right] \leq \tilde C_n e^{C_n |x|}\E  \left[e^{C_n W_1^{*}}\right] <\infty
 $$
for some constant $\tilde C_n  >  0$. Thus, for every $x$ the derivative $(\partial^n/\partial y^n) Q^{(0)} (x, B^{-1}_{x}(y) )$ admits an  integrable bound not depending on $y$. This  implies that the function
$y\mapsto \rho_x(y)$ is  infinitely differentiable and 
$$
\frac{\partial^n}{\partial y^n}\rho_x(y)=\E \left[\frac{\partial^n}{\partial y^n} Q^{(0)} (x, B^{-1}_{x}(y) )\right].
$$
Moreover, increasing in the need the constant $C_n$ we obtain  the inequality 
$$
\sup_y \left |\frac{\partial^n}{\partial y^n}\rho_x(y)\right |\le  C_n e^{C_n |x|}.
$$

Since $B^{-1}_x(0+) = \infty$, the inequality (\ref{Qn}) implies that $Q^{(n)}(x, B^{-1}_x(0+)) = 0$. Hence, $(\partial^n / \partial y^n) \rho_x(0) = 0$. 

Differentiating the function $u\mapsto J(u,t)$ we have: 
$$ 
\frac{\partial^{n}}{\partial u^n} J(u,t)=  \int \left[\int  H(t,x,y)\frac{\partial^{n}}{\partial u^n}\rho_x(y-u)dy\right]\varphi_{0,1}(x)dx. 
$$

\subsection{Smoothness of $\Upsilon^1$.}

Now we get a result on smoothness  of the function 
$$
\Upsilon^{1}(u)= \int_0^2  \int_{D}  \E\left[\Phi(Y_t^{u}(z)+\xi)]\right]
m(t,dz) \alpha(dt).
$$
Recall that $c>0$ and $\xi<0$.  For us  it is more convenient to work with the random variable $\zeta:=-\xi>0$. The needed result 
follows from

\begin{proposition}
Let $Y^u=e^{V}(u+c e^{-V}\cdot \Lambda)$  where $V=\kappa\cdot \Lambda+\sigma W$,  $\kappa$ and $\sigma>0$ are measurable functions.   
 Let $\zeta>0$ be a random variable with a density $f$  two times  continuously  differentiable on $ (0, \infty)$  and such that   $f',f''\in L^1(\R_{+})$.
Then the function $\phi_t(u):=\E\left[\Phi(Y_t^{u}- \zeta)]\right]$   is two times  continuously differentiable in $u$ and there is a constant $C$ depending only of $\kappa^*$ and $\sigma^*$ and  such that 
 \beq
 \label{bounds}
 |\phi_t'(u)|  \le C ||f'||_{L^1}, \qquad |\phi_t''(u)|\le C ||f''||_{L^1}.
 \eeq
 for all $t\in [0,2]$, $u>0$. 
\end{proposition} 

\noindent
{\sl Proof.} 
Let $H:\R\to [0,1]$ be a measurable function such that $H= 0$ on $\R_-$ and let $\gamma>0$ be a constant. Then 
the function 
$$
h(y):= \int_{\R_+} H(y-\gamma w)f(w)dw =(1/\gamma)\int_{\R_+} H(z)f((y-z)/\gamma )dz. 
$$
It follows that $h$ is two times differentiable,   
$$
h'(y)=(1/\gamma^2)\int_\R H(z)f'((y-z)/\gamma)dz, \qquad  h''(y)=(1/\gamma^{ 3 })\int_\R H(z)f''((y-z)/\gamma)dz, 
$$ 
$|h'(y)|\le  (1/\gamma) ||f'||_{L^1}$ and $|h''(y)|\le (1/\gamma^2) ||f''||_{L^1}$. 

Let us transform our problem to the above elementary framework. 
Using the definitions of $Y^u$ and $V$
we get that 
$$
\phi_t(u)=\E [\Phi (Y^u_t-\zeta)]=\E\big [
 \Phi_t\big (e^{\sigma\cdot W_t} {\color{red} (}  u-e^{-\sigma\cdot W}m\cdot \Lambda _t {\color{red} )} - e^{-\kappa\cdot \Lambda_t}\zeta \big )\big ],
$$
where the functions $\Phi_t(y):=\Phi(e^{ \kappa\cdot \Lambda_t}y)$, 
$m:=c e^{-\kappa\cdot \Lambda_t}$. Thus, 
$$
\phi_t(u) =\int_{\R} \E\big [\Phi_t\big (e^{\sigma\cdot W_t} {\color{red} (} u-e^{-\sigma\cdot W}m\cdot \Lambda _t {\color{red} )} -e^{-\kappa\cdot \Lambda_t}w\big )\big ]   f(w)dw.
$$

To get rid of stochastic integrals in the above formula we consider on $\R_+$ the strictly increasing function $A:=\sigma^2\cdot \Lambda $  with the inverse $A^{-1}$ and observe that the law of the process 
$\sigma\cdot W$ is the same as of the process $W_A$. Changing the variable  in the integral with respect to $\Lambda$ we get that 
$e^{W_A}m\cdot \Lambda_{t}=e^{W}\tilde m\cdot \Lambda_{A_t}$ where the function $\tilde m(r):=m(A^{-1}_r)/\sigma^
2_{A^{-1}_{r}}$. 

It follows that  
$$
\E\big [\Phi_t\big (e^{\sigma\cdot W_t}  (u-e^{-\sigma\cdot W}m\cdot \Lambda_t  )  -e^{-\kappa\cdot \Lambda_t}w\big )\big ] =\E\big [\Phi_t\big (e^{W_{A_t}}  ( u-e^{ -W}\tilde  m \cdot \Lambda_{A_t} )  -e^{-\kappa\cdot \Lambda_t}w\big )\big ].  
$$  
The change of variable $s\to st/A_t$  replace the integration over the interval $[0,A_t]$ by the integration over the interval $[0,t]$ and the rhs of the above equality is equal to 
\beq
\label{rhs}
\E\big [\Phi_t\big (e^{\sqrt {A_t /t} W_{t}}  ( u-e^{-\sqrt {A_t/t}W}\bar  m_t \cdot \Lambda_{t}   ) -e^{-\kappa\cdot \Lambda_t}w\big )\big ]
\eeq
where 
$\bar m_t(s):=\tilde m(A_ts)A_t = m(T_t(s))A_t/\sigma^2_{T_t(s)}$ with the abbreviatio  $T_t(s):=A^{-1}_{A_ts}$.  

Note that in the conventional notation 
$$
e^{-\sqrt {A_t/t}W}\bar  m_t \cdot \Lambda_{t}:=\int_0^t e^{-\sqrt {A_t/t}W_s}\bar m_t(s)ds.
$$

Using the conditioning with respect to $W_t=x$  and taking into account that the conditional law of $(W_s)_{s\in[0,t]}$ is that of the Brownian bridge on $[0,t]$ ending at the level $x$, that is the same as of the process
$W_s+(s/t)(x-W_t)$, $s\in [0,t]$,  
 we infer that (\ref{rhs}) can be written as 
 $$
\frac 1{\sqrt {t}} \int_{\R}  H (u-\gamma_{t,x} w)\varphi_{0,1} (x/\sqrt t)dx
 $$ 
 where $H(y)=H_{t,x}(y):=\E [\Phi_t(e^{\sqrt {A_t/t}x}(y- \eta_{t,x})]$, $\gamma_{t,x}:=e^{-\kappa\cdot \Lambda_t-\sqrt {A_t/t}x}$, 
 $$
  \eta_{t,x}:=\int_0^t e^{-\sqrt {A_t/t}(W_s-(s/t)W_t)}m_{t,x}(s)ds, \quad m_{t,x}(s):=e^{(s/t)x-\sqrt {A_t/t}x}\bar m_t(s).
 $$
Thus, 
$$
\phi_t(u)=\frac 1{\sqrt {t}} \int_{\R} \left[\int_{\R}  H (u-\gamma_{t,x} w)f(w)dw\right] \varphi_{0,1} (x/\sqrt t)dx.
$$

Let denote by $h(u)=h_{t,x}(u):=\E[ H (u-\gamma_{t,x} \zeta)]$, i.e. the function in parentheses above. Then $h$ is two times continuously differentiable  and  there is a constant $C$ such that for all $t\in [0,2]$ and all $x$  
\bean
|h'(u)|&\le &  (1/\gamma_{t,x}) ||f'||_{L^1}\le C e^{C|x|} ||f'||_{L^1},\\
 |h''(u)|&\le&  (1/\gamma^2_{t,x}) ||f''||_{L^1}\le   C e^{C|x|} ||f''||_{L^1}.
\eean
Since 
$$
\frac 1{\sqrt {t}} \int_{\R}  e^{C|x|} \varphi_{0,1} (x/\sqrt t)dx=\int_{\R} e^{2C\sqrt{t} |x| } \varphi_{0,1} (x)dx  \le \int_{\R}   e^{4C |x| }\varphi_{0,1} (x)dx < \infty,
$$
this implies that the function $\phi_t(u)$ is two times continuously differentiable in $u$ and  the bounds  (\ref{bounds}) 
hold. 
\smallskip

\begin{remark} Minor changes in the above proof allows to get smoothness result for annuity model and for the model where the distribution $F_\xi$ charges both half-axes. 
\end{remark}

\section{Integro-differential equations}
Knowing that the survival  probability is a $C^2$-function, the derivation of the integro-differential equation is easy and 
we get it for all variants of the model.  
\begin{proposition}
\label{propos1}
Assume that $\Phi=(\Phi_0,\Phi_1,\dots,\Phi_{K-1})\in C^2$. Then  $\Phi$
satisfies the system of integro-differential equations:
\beq
\label{IDE}
 \frac 12 \sigma^2_i u^2 \Phi''_i(u)+ (a_i u+c)\Phi'_i(u) + \mathcal{I}(\Phi_i)(u) + \sum_{j\in E}\lambda_{i j}\Phi_j(u), \quad i\in E,
\eeq
where 
\begin{eqnarray*}
\mathcal{I}(\Phi_i)(u) = \alpha\int (\Phi_i(u +x) - \Phi_i(u))F_{\xi}(dx)= \alpha\int \Phi_i(u +x) F_{\xi}(dx)- \alpha \Phi_i(u).
\end{eqnarray*}
\end{proposition}
\noindent{\sl Proof.}
Take arbitrary $h > 0$ and $\e>0$ such that $u\in [\e,\e^{-1}]\subset (0,1)$.  
Skipping, as usual, $i$ in the notations we define the stopping time 
$$
\tau^{\e}_h:=\inf\big\{t\ge 0\colon\  X_{t}
\notin[\e,,\e^{-1}]\big\}\wedge T_1\wedge \tau_1\wedge h. 
$$

The  Ito formula applied to the function $\Phi$ and the semimartingale $(X,\theta) = (X^{u, i},\theta^i)$ yields in the identity 
\beq
\label{Id}
 \Phi_{\theta_t} (X_t)
 =
 \Phi_i(u)+\Phi'_{\theta_-}(X_-)\sigma \cdot W_t +\cL \Phi_\theta(X)\cdot \Lambda_t+\sum_{s\le t} \big(\Phi_{\theta_s}(X_s)-\Phi_{\theta_{s-}}(X_{s-})\big )
\eeq
where 
$$
\cL \Phi_\theta(X):=\frac 12\sigma_\theta^2
  X^2\Phi''_{\theta}(X)+  (a_{\theta}X+c)\Phi'_{\theta}(X),  \quad \Sigma_t:=\sum_{s\le t} \big(\Phi_{\theta_s}(X_s)-\Phi_{\theta_{s-}}(X_{s-})\big ).
$$
 Due to independence, the processes $P$ and $\theta$ have no common jumps and 
\bean
 \Sigma_t
 &=&\int_0^t\int_{\R} (\Phi_{\theta_{s-}}(X_{s-} + x)-\Phi_{\theta_{s-}}(X_s-))\pi(ds,dx)\\
 &&+\int_0^t\sum_{k\in E} (\Phi_{k}(X_{s-})-\Phi_{\theta_{s-}}(X_s-))p(ds,k).
 \eean
 
 Now we replace in the formula (\ref{Id}) $t$ by $\tau^{\e}_h$ and take the expectation from both sides using the following observations: 
 
-- In virtue of  the strong Markov property $\Phi_i(u)= {\bf E} \Phi_{\tau^{\e}_h}(X_{\tau^{\e}_h})$ (recall that $\Phi_i(u)=0$ for $u\le 0$). 

 -- For any $\e>0$, the integrands on $[0,\tau^{\e}_h(\omega)]$ are bounded from above  and, therefore, the expectation of the stochastic integral over the Wiener process is zero. 

--  As ${\tau^{\e}_h}\le \tau_1$, then we can replace $\theta_{s-}$ in the integrals by $i$. In addition, $\tau^\e_h=h$ when $h$ is small enough.  

-- According to the definition of dual predictable projections
$$ 
\E \int_{0}^{\tau^{\e}_h}\int(\Phi_i(X_{s-}+x) - \Phi_i(X_{s-})) \pi(ds,dx) = \E \int_{0}^{\tau^{\e}_h}\int(\Phi_i(X_{s-}+x) - \Phi_i(X_{s-}))ds \Pi(dx),
$$
$$
 \E \int_{0}^{\tau^{\e}_h}\sum_{j\in E\setminus \{i\}}\Phi_j(X_{s-}-\Phi_i(X_{s-}))  p(ds,j) = 
 \E \int_{0}^{\tau^{\e}_h}\sum_{j\in E}\Phi_j(X_{s-}))\lambda _{ij}ds.  
$$
As a result we get that 
$$
\E \int_{0}^{\tau^{\e}_h}\Bigg( \cL \Phi_i(X_s)
   +\int(\Phi_i(X_{s}+x) - \Phi_i(X_{s}))F_\xi(dx)+\sum_{j\in E}(\Phi_j(X_{s})- \Phi_{i}(X_{s}))\lambda _{ij}\Bigg)ds=0.
$$

Dividing both sides of this equality by $h$ and letting $h\downarrow 0$ we get the result. \fdem  
\begin{remark}
The obtained equation is homogeneous and holds as well for the ruin probability $\Psi(u)=1-\Phi(u)$. 
\end{remark}

\section{Exponential distribution and differential equations}
If  $F_\xi(dx)=\mu e^{-x/\mu}dx$, then the ruin probability $\Psi\in C^2$. Note that 
$$
\frac {d}{du}\int_0^\infty \Psi_i(u+y)e^{-y/\mu}dy=-\Psi_i(u)+\frac 1\mu \int_0^\infty \Psi_i(u+y)e^{-y/\mu}dy
$$
and therefore, 
$$
\frac {d}{du}\cI\Psi_i(u)=\frac 1\mu \cI\Psi_i(u)-\alpha \Psi'_i(u). 
$$
 
Differentiating the integro-differential equation   for $\Psi_i$
(it is  same as for $\Phi_i$) and subtracting from the resulting identity multiplied by $\mu$ the equation
(\ref{IDE})
obtain an ODE of the third order which can be transformed to the form  
$$
\label{stoch_sys}
\Psi'''_i- p_i(u)\Psi''_i+ q_i(u)\Psi'_i + \frac {2}{ \mu \sigma_i^2}\frac 1{u^2} \sum_{j\in E}\lambda_{i j}\Psi'_ j(u)  - \frac {2}{\mu\sigma_i^2}\frac 1{u^2}  \sum_{j\in E}\lambda_{i j}\Psi_j(u)=0,
$$
where 
$$
p_i(u):= \frac {1}{\mu}-2\Big(1+\frac{a_i}{\sigma_i^2}\Big)\frac {1}{u}-\frac{2c}{\sigma_i^2}\frac {1}{u^2},\qquad 
q_i(u):=-\frac {2a_i}{\mu \sigma_i^2} \frac {1}{u}+(a_i-\alpha -c/\mu)\frac {2}{\sigma_i^2}\frac {1}{u^2}.
$$
\smallskip
{\bf Acknowledgement.} This work was supported by the Russian Science Foundation associated grants  20-68-47030 and 20-61-47043.  

\smallskip
The authors express their thanks to Paul Embrechts who encouraged them  to do a research on the smoothness of the ruin probabilities.


\begin{thebibliography}{100}
\bibitem{DKR} 
Di Masi, G.B., Kabanov, Yu.M., Runggaldier, W.J.,  
Mean-square hedging of options on a stock with Markov volatilities.  {\em Theory Probab. Appl.} {\bf 39} (1) (1994), 172--182.  


\bibitem{EK} 
Ellanskaya, A., Kabanov, Yu., On ruin probabilities with risky investments in a stock with stochastic volatility. {\em Extremes} 
{\bf 24} (4) (2021) 687--697. 


\bibitem{KP2016}
Kabanov, Yu., Pergamenshchikov, S.,  In the  insurance business risky investments are dangerous: the case of negative risk sums. {\em Finance and Stochastics},   {\bf 20} (2) (2016),  355--379. 

\bibitem{KP2022} 
Kabanov, Yu., Pergamenshchikov, S.,
On ruin probabilities with investments in a risky asset with a regime-switching price.
{\em Finance and Stochastics} {\bf 26} (4) (2022), 877--897. 

\bibitem{KPukh} 
Kabanov, Yu., Pukhlyakov, N., Ruin probabilities with investments:
smoothness, IDE and ODE, asymptotic behavior.  {\sl J. Appl. Probab.}  {\bf 59} (2) (2020), 556--570. 

\bibitem{Paul-93} 
Paulsen, J.,
 Risk theory in stochastic economic environment.  
 {\em Stoch. Proc. Appl.} 
 {\bf 46} (2) (1993), 327--361. 
  
\bibitem {RP}
Ramsden, L., Papaioannou, A.,
Asymptotic results for a Markov-modulated risk process with stochastic investment. {\em J. Comp. Appl. Math.} {\bf 313} (2017), 38--53. 
\end{thebibliography}
\end{document}